\documentclass[10pt]{article}
\usepackage{amsmath,amsthm}
\usepackage{latexsym}
\usepackage{amssymb,amscd}
\usepackage{xspace}
\usepackage{amscd}
\usepackage{enumerate}
\newtheorem{theo}{Theorem}[section]
\newtheorem{conj}[theo]{Conjecture}

\newtheorem{rema}[theo]{Remark} 

\newtheorem{lemm}[theo]{Lemma}
\newtheorem{coro}[theo]{Corollary}

\newcommand{\ch}{\S}

\newcommand{\sumind}[2]{\genfrac{}{}{0pt}{}{#1}{#2}}
\numberwithin{equation}{section}
\title{Towards large genus asymtotics of intersection numbers on moduli spaces of curves}
\author{Maryam Mirzakhani, Peter Zograf} 
\begin{document}
\maketitle
\begin{abstract}
We explicitly compute the diverging factor in the large genus asymptotics of the Weil-Petersson volumes of the moduli spaces 
of $n$-pointed complex algebraic curves. Modulo a universal multiplicative constant we prove the existence of a complete asymptotic expansion of the Weil-Petersson volumes in the inverse powers of the genus with coefficients that are polynomials
in $n$. This is done by analyzing various recursions for the more general intersection numbers of tautological classes, whose large genus asymptotic behavior is also extensively studied.
\end{abstract}

\begin{section}{Introduction and statement of results}
In this note, we study the asymptotic behavior of the Weil-Petersson volumes $V_{g,n}$ of the moduli 
spaces $\mathcal{M}_{g,n}$ of $n$-pointed complex algebraic curves of genus $g$ as $g \rightarrow \infty.$ Here
$$V_{g,n}=\frac{1}{(3g-3+n)!} \int_{{\mathcal{M}}_{g,n}} \omega_{g,n}^{3g-3+n},$$ where $\omega_{g,n}$ is the Weil-Petersson symplectic form on $\mathcal{M}_{g,n}.$ 

The following conjecture was made in \cite{Z:con} on the basis of numerical data:

\begin{conj}\label{conj}
For any fixed $n\geq 0$ 
$$V_{g,n} =(2g-3+n)! (4\pi^2)^{2g-3+n} \frac{1}{\sqrt{g \pi}} \left(1 + \frac{c_{n}}{g} + O\left(\frac{1}{g^{2}}\right)\right)$$
as $g\rightarrow \infty.$ 
\end{conj}
 
The objective of this paper is to prove the statements formulated below.

\begin{theo}\label{theo:main:fixn}
 There exists a universal constant $C\in (0,\infty)$ such that for any given $k \geq 1, n\geq 0,$ 
 $$V_{g,n} =C\,\frac{(2g-3+n)!\,(4\pi^2)^{2g-3+n}}{\sqrt{g}} \left(1+\frac{c_n^{(1)}}{g}+\ldots+\frac{c_n^{(k)}}{g^{k}}+ O\left(\frac{1}{g^{k+1}}\right)\right), $$
 as $g \rightarrow \infty.$
Each term $c_n^{(i)}$ in the asymptotic expansion is a polynomial in $n$ of degree $2i$ with coefficients in 
${\mathbb Q}[\pi^{-2},\pi^2]$ that are effectively computable.  Moreover, the leading term of $c_n^{(i)}$ is equal to
$\frac{(-1)^i}{i!\,(2\pi^2)^i}\,n^{2i}.$
\end{theo}

\begin{rema}\label{rem1}
{\rm Note that Conjecture $\ref{conj}$ claims that $C=\frac{1}{\sqrt{\pi}}$. 
Numerical data suggest that the coefficients of $c_n^{(i)}$ actually belong to ${\mathbb Q}[\pi^{-2}]$. For example,}
$$c_n^{(1)}=-\frac{n^2}{2\pi^2}-\left(\frac{1}{4}-\frac{5}{2\pi^2}\right)n+\frac{7}{12}-\frac{17}{6\pi^2}.$$ 
\end{rema}

Our method also implies that given $k\geq 1, n \geq 0$ we have
\begin{equation}\label{1}
\frac{V_{g,n+1}} {8\pi^2 gV_{g,n}}= 1+ \frac{a_n^{(1)}}{g}+\ldots+\frac{a_n^{(k)}}{g^{k}}+O\left(\frac{1}{g^{k+1}}\right),
\end{equation}
\begin{equation}\label{2}
\frac{V_{g-1,n+2}}{V_{g,n}}=1+ \frac{b_n^{(1)}}{g}+\ldots+\frac{b_n^{(k)}}{g^{k}}+O\left(\frac{1}{g^{k+1}}\right),
\end{equation}
as $g\to\infty$, where the coefficients $a_n^{(i)}$ and $b_n^{(i)}$ can be explicitly computed. However, here we do this only for $a_n^{(1)}$ and $b_n^{(1)}$:
\begin{theo}\label{11}
For any fixed $n \geq 0$: 
 \begin{align}\label{110}
 &\frac{V_{g,n+1}}{8\pi^2 gV_{g,n}}= 1+\left(\left(\frac{1}{2}-\frac{1}{\pi^2}\right)n-\frac{5}{4}+\frac{2}{\pi^2}\right)\cdot\frac{1}{g}+ O\left(\frac{1}{g^2}\right),\\
 \label{111}
 &\frac{V_{g-1,n+2}}{V_{g,n}}=1+\frac{3-2n}{\pi^2}\cdot\frac{1}{g}+O\left(\frac{1}{g^2}\right).
 \end{align}
\end{theo}

With the help of the identity $\frac{V_{g+1,n}}{V_{g,n}}=\frac{V_{g+1,n}}{V_{g,n+2}}\cdot\frac{V_{g,n+2}}{V_{g,n+1}}\cdot\frac{V_{g,n+1}}{V_{g,n}}$ this theorem immediately yields
\begin{coro}\label{c}
Let $n \geq 0$ be fixed, then
$$\frac{V_{g+1,n}}{V_{g,n}}= (4\pi^2)^2(2g+n-1)(2g+n-2)\left(1-\frac{1}{2g}\right)\cdot \left(1+ \frac{r_{n}^{(2)}}{g^2}+ O\left(\frac{1}{g^3}\right)\right),$$
as $g \rightarrow \infty.$
\end{coro}

Note that  $r_n^{(2)}=3/8-c_n^{(1)}$ by Theorem \ref{theo:main:fixn}, so that
$$r_n^{(2)}=\frac{n^2}{2\pi^2}+\left(\frac{1}{4}-\frac{5}{2\pi^2}\right)n-\frac{5}{24}+\frac{17}{6\pi^2}.$$

\begin{rema}\label{rem2}
{\rm Since $\prod_{g=1}^{\infty} (1+a_{g})$
converges when $a_{g}=O(1/g^2)$, Corollary \ref{c} easily implies that there exists $$\lim_{g\rightarrow \infty} \frac{V_{g,n} \sqrt{g}}{(2g-3+n)!\,(4\pi^2)^{2g-3+n}} =C\in (0,\infty).$$
The proof of Theorem $\ref{theo:main:fixn}$ is based on $(\ref{1})$ and $(\ref{2})$ (see Theorems $\ref{theo:general:1}$,
$\ref{theo:general:2}$) and follows similar lines as well. The method we use also allows to calculate the error terms explicitly: even though we only do the calculation for the coefficient of $1/g$, the error terms of order of $1/g^s$ can be written in terms of values of the intersections of $\psi$ classes on surfaces of genus at most $s$.  However, this method does not provide any information about the exact value of $C$.}
\end{rema}

Analyzing the signs of the error terms of order $1/g^2$ in Theorem $\ref{11}$, we get
\begin{coro}\label{m}
Given $n \geq 2,$ there exists $g_{0}$ such that the sequence 
 $$\left \{\frac{V_{g-1,n+2}} {V_{g,n}}\right\}_{g \geq g_{0}}$$
 is increasing. 
 Similarly, for $n \geq 3$ there exists $g_{0}$ such that the sequence 
$$ \left\{ \frac{8\pi^2 g V_{g,n}} {V_{g,n+1}}\right\} _{g \geq g_{0}}$$
is increasing.
\end{coro}

We also obtain somewhat weaker results when $n$ varies as $g\to \infty$:

\begin{theo}\label{Main}
For any sequence $\{n(g)\}_{g=1}^{\infty}$ of non-negative integers with 
$$\lim_{g\rightarrow \infty} \frac{n(g)^2}{g}=0,$$ 
we have
$$V_{g,n(g)}=  \frac{C}{\sqrt{g}}\, (2g-3+n(g))!\,(4\pi^2)^{2g-3+n(g)}\, \left(1 + O\left(\frac{1+n(g)^2}{g}\right)\right)$$
as $ g\rightarrow \infty.$ 
\end{theo}
In fact, we prove that 
$$ \frac{\sqrt{g} \,V_{g,n(g)}}{C\,(2g-3+n(g))!\;(4\pi^2)^{2g-3+n(g)}}= O\left(\frac{(g+n)!}{g! \, g^{n}}\right). $$
 
\bigskip
\noindent
{\bf Notes and remarks.}
\begin{itemize}
\item 
\noindent
It may be instructive to compare Theorem \ref{theo:main:fixn} to 
the asymptotic formula for the Weil-Petersson volumes for fixed $g\geq 0$ and $n\to\infty$ (cf. \cite{MZ}, Theorem 6.1):

\begin{theo}\label{z}
 For any fixed $g\geq 0$ 
\begin{equation}\label{zn} 
V_{g,n} = n! C^{n} n^{(5g-7)/2} \left(c_{g}^{(0)} + \frac{c_{g}^{(1)}}{n} + \dots\right),
\end{equation}
as $n \rightarrow \infty$, where $C=-\frac{2}{x_0\,J'_0(x_0)}$, $x_0$ is the first positive zero of the Bessel function $J_0$, and the coefficients $c_{g}^{(0)},c_{g}^{(1)},\dots$ are effectively computable.
\end{theo}
 \item Penner \cite{P:vol} developed a different method for computing Weil-Petersson volumes by integrating the Weil-Petersson volume form over simplices in the cellular decomposition of the moduli space. In \cite{Gr:vol}, this method of integration was used to prove that for a fixed $n>0$ there are constants $C_{1},C_{2}>0$ such that 
 $$C_{1}^g \cdot (2g)! < V_{g,n} < C_{2}^g\cdot (2g)!.$$ 
(This result was extended to the case of $n=0$ by an algebro-geometric argument in \cite{ST:vol}.) 
 Note that these estimates do not give much information about the growth of $V_{g-1,n+2}/V_{g, n}$ and $V_{g,n}/V_{g,n+1}$ when 
 $g \rightarrow \infty.$
\item The estimates from $\cite{M:large}$ imply that given $n \geq 0$ there exists $m>0$ such that 
\begin{equation}\label{compare}
 g^{-m} \leq \frac{V_{g,n}}{(2g-3+n)! \,(4\pi^2)^{2g-3+n)} }\leq g^{m}. 
\end{equation}
In general, $$\lim _{g+n \rightarrow \infty} \frac{\log(V_{g,n})}{(2g+n)\log(2g+n)}=1,$$
but understanding the asymptotics of $V_{g,n}$ for arbitrary $g,n$ seems to be more complicated.
\end{itemize}
\end{section}

\begin{section}{Relations between intersection numbers}\label{Asym} 

To begin with, let us recall some well-known facts about tautological classes on $\overline{\mathcal{M}}_{g,n}$ and their intersections. For $n>0,$ there are $n$ tautological line bundles $\mathcal{L}_{i}$ on $\overline{\mathcal{M}}_{g,n}$ whose fiber at the point 
$(C,x_{1},\ldots,x_{n})\in \overline{\mathcal{M}}_{g,n}$ is the cotangent line to $C$ at $x_{i},$ and we
put $\psi_{i}= c_{1}(\mathcal{L}_{i}) \in H_{2} (\overline{\mathcal{M}}_{g,n}, {\mathbb Q})$ (cf. e.g. \cite{Harris:book} or \cite{AC}).

\noindent
{\bf Notation.} For $d=(d_{1},\ldots,d_{n})$ with $ d_{i} \in {\mathbb Z} _{\geq 0}$ put $|d|=d_{1}+\ldots+d_{n}$ and, assuming $|d| \leq 3g-3+n,$ put $d_{0}= 3g-3+n-|d|$. Define
\begin{align*}
[\tau_{d_{1}}\ldots \tau_{d_{n}}]_{g,n} =\left[ \prod_{i=1}^{n} \tau_{d_{i}}\right]_{g,n}&= \;\frac{\prod_{i=1}^{n} 2^{2d_i}\,(2d_{i}+1)!!} {d_{0}!} \int_{\overline{\mathcal{M}}_{g,n}} \psi_{1}^{d_{1}} \cdots \psi_{n}^{d_{n}} \omega_{g,n}^{d_{0}}\\
&=\prod_{i=1}^{n} 2^{2d_i}\,(2d_{i}+1)!!\; \frac{(2\pi^2)^{d_0}}{d_0!} \int_{\overline{\mathcal{M}}_{g,n}} \psi_{1}^{d_{1}} \cdots \psi_{n}^{d_{n}}\kappa_1^{d_0},
\end{align*}
 where $\kappa_1=\frac{[\omega_{g,n}]}{2\pi^2}$ is the first Mumford class on $\overline{\mathcal{M}}_{g,n}$. 
According to $\cite{M:JAMS},$ for $L=(L_{1},\ldots,L_{n})$ the Weil-Petersson volume of the moduli space 
of hyperbolic surfaces of genus $g$ with $n>0$ geodesic boundary components of lengths $2L_{1},\ldots ,2L_{n}$ can be written as
\begin{equation}\label{re2}
V_{g,n}(2 L)=\sum_{\sumind{d_1,\dots, d_n}{|d| \leq 3g-3+n}} \left[ \prod_{i=1}^{n} \tau_{d_{i}}\right]_{g,n} \frac{L_{1}^{2d_{1}}}{(2d_{1}+1)!}\cdots \frac{L_{n}^{2 d_{n}}}{(2d_{n}+1)!}.
\end{equation} 

\noindent
{\bf Recursive formulas.} 
The following recursions for the intersection numbers  $[\tau_{d_{1}}\ldots \tau_{d_{n}}]_{g,n}$ hold:

\bigskip
\noindent

\begin{align*}{\rm ({\bf Ia})}\qquad\quad\;
 \left[ \tau_{0} \tau_{1} \prod_{i=1}^{n} \tau_{d_{i}}\right]_{g,n+2}=\quad &\left[\tau_{0}^{4} \prod_{i=1}^{n} \tau_{d_{i}}\right]_{g-1,n+4} +\\
+\; 6 \sum_{\genfrac{}{}{0pt}{}{g_{1}+g_{2}=g}{I\amalg J=\{1,\ldots, n\}}} &\left[ \tau_{0}^{2} \; \prod_{i\in I } \tau_{d_{i}}\right]_{g_{1},|I|+2} \cdot \quad\left[\tau_{0}^{2}\prod_{i\in J} \tau_{d_{i}}\right]_{g_{2},|J|+2}.
\end{align*}

\noindent
\begin{align*}{\rm ({\bf Ib})}\qquad\quad
\left[\tau_{0}^2 \tau_{l+1} \prod_{i=1}^{n} \tau_{d_{i}}\right]_{g,n+3}
&=\quad \left[\tau_{0}^{4} \tau_{l} \prod_{i=1}^{n} \tau_{d_{i}}\right]_{g-1,n+5} +\\
+\;8 \sum_{\genfrac{}{}{0pt}{}{g_{1}+g_{2}=g}{I\amalg J=\{1,\ldots, n\}}} 
&\left[ \tau_{0}^{2} \tau_{l} \; \prod_{i\in I } \tau_{d_{i}}\right]_{g_{1},|I|+3} 
\cdot \quad\left[\tau_{0}^{2}\prod_{i\in J} \tau_{d_{i}}\right]_{g_{2},|J|+2}+\\
+\; 4\sum_{\genfrac{}{}{0pt}{}{g_{1}+g_{2}=g}{I\amalg J=\{1,\ldots, n\}}}
&\left[ \tau_{0}\tau_{l} \; \prod_{i\in I } \tau_{d_{i}}\right]_{g_{1},|I|+2} 
\cdot\quad\left[\tau_{0}^{3}\prod_{i\in J} \tau_{d_{i}}\right]_{g_{2},|J|+3}.
\end{align*}

\medskip
\noindent
$$
{\rm ({\bf II})}\quad (2g-2+n) \left[ \prod_{i=1}^{n} \tau_{d_{i}}\right]_{g,n}= \quad\frac{1}{2} \sum_{l=1}^{3g-2+n} \frac{(-1)^{l-1}\, l\,\pi^{2l-2}}{(2l+1)!} \left[ \tau_{l} \; \prod_{i=1}^{n} \tau_{d_{i}}\right]_{g,n+1}.\quad
$$

\medskip
\noindent
({\bf III})
Put $a_{i} =(1-2^{1-2i})\,\zeta(2i),$ where $\zeta$ is the Riemann zeta function and $i\in\mathbb{Z}_{\geq 0}$.
Then
$$
[ \tau_{d_{1}}\ldots\tau_{d_{n}}]_{g,n}= {A}_{{d}} + {B}_{{d}}+ {C}_{{d}},
$$
where 
\begin{align}
\label{A}
{A}_{{d}}=&8 \; \sum_{j=2}^{n} \sum_{l=0}^{d_{0}} (2d_{j}+1) \; a_{l} \left[\tau_{d_{1}+d_{j}+l-1} \prod_{i\not=1,j} \tau_{d_{i}}\right]_{g,n-1},\\
\label{B}
{B}_{{d}}= &16 \;\sum_{l=0}^{d_{0}} \sum_{\genfrac{}{}{0pt}{}{k_{1}+k_{2}=}{=l+d_{1}-2}} a_{l} \left[\tau_{k_{1}} \tau_{k_{2}} \prod_{i\not=1} \tau_{d_{i}}\right]_{g-1,n+1},\\
\label{C}
{C}_{{d}}=&16\sum_{\genfrac{}{}{0pt}{}{g_{1}+g_{2}=g}{I\amalg J=\{1,\ldots, n\}}} \sum_{l=0}^{d_{0}} \sum_{\genfrac{}{}{0pt}{}{k_{1}+k_{2}=}{=l+d_{1}-2}} a_{l} \; \left[\tau_{k_{1}} \prod_{i\in I } \tau_{d_{i}}\right]_{g_1,|I|+1} \cdot \left[\tau_{k_{2}} \prod_{i\in J} \tau_{d_{i}}\right]_{g_2,|J|+1}.
\end{align}

\noindent
{\bf Basic properties of the sequence ${\{a_{i}=(1-2^{1-2i}) \zeta(2i)\}.}$}
It is easy to check that for $i \geq 1$
$$a_i =\frac{1}{(2i-1)!} \int_{0}^{\infty} \frac{t^{2i-1}}{1+e^{t}}\; dt$$
and
\begin{equation}\label{adif}
a_{i+1}-a_{i}=\int_{0}^{\infty} \frac{1}{(1+e^{t})^{2}} \left(\frac{t^{2i+1}}{(2i+1)!}+\frac{t^{2i}}{2i!}\right) dt.
\end{equation}
\begin{lemm}\label{aa}
The sequence $\{a_{i}\}_{i=1}^{\infty}$ is increasing. 
Moreover,
\begin{enumerate}[(i)]
\item 
$$\lim_{i\rightarrow \infty} a_{i}=1,\qquad \sum_{i=0}^{\infty} (a_{i+1}-a_{i})=\frac{1}{2},$$
\item 
$a_{i+1}-a_i$ has the order of $1/2^{2i}$, i.e., there exist $C>0$ such that
\begin{equation}\label{abound}
\frac{1}{C\cdot 2^{2i}}<a_{i+1}-a_{i} <\frac{C}{2^{2i}},
\end{equation}
\item 
\begin{equation}\label{asums}
\sum_{i=0}^{\infty} i (a_{i+1}-a_{i})=\frac{1}{4},
\end{equation}
\item for $j \in {\mathbb Z}, j \geq 2,$ the sum
$$\sum_{i=0}^{\infty} i^{j} (a_{i+1}-a_{i})$$
is a polynomial in $\pi^2$ of degree $[j/2]$ with rational coefficients.
\end{enumerate}
\end{lemm}

\noindent
{\bf Proof.}
Both $(i)$ and $(ii)$ easily follow from the definition of $a_{i}$ and $(\ref{adif})$. 
As for $(iii)$, let 
$$S_{1}=\sum_{i=0}^{\infty} \int_{0}^{\infty} \frac{1}{(1+e^{t})^{2}}\cdot \frac{t^{2i+1}}{(2i+1)!}\,dt$$
and
$$S_{2}= \int_{0}^{\infty} \frac{t \; e^t}{(1+e^{t})^{2}} dt.$$
We have 
$$S_{1}=  \int_{0}^{\infty} \frac{(e^{t}-e^{-t})\,dt}{2(1+e^{t})^{2}}= -\frac{1}{2}+\log 2\,,\qquad S_{2}=\log 2,$$
so that from $(\ref{adif})$, $$2 \sum_{i=1}^{\infty} i(a_{i+1}-a_{i})+S_{1}=S_{2},$$ which implies $(\ref{asums}).$

We will prove $(iv)$ by induction in $j$. The base case $j=1$ being checked in $(iii)$, we observe that
\begin{align*} 
\sum_{i=0}^{\infty} \frac{i^{j}\; t^{2i+1}}{(2i+1)!}&=(t D t^{-1})^j(\sinh t) 
=\sum_{l=0}^j t^l(a_{1,j}^{(l)}\cosh t + b_{1,j}^{(l)}\sinh t),\\
\sum_{i=0}^{\infty} \frac{i^{j} t^{2i}}{(2i)!}&=D^j (\cosh t)=\sum_{l=0}^j t^l(a_{2,j}^{(l)}\cosh t + b_{2,j}^{(l)}\sinh t),
\end{align*}
where $D=\frac{t}{2}\cdot\frac{d}{dt}$, and the coefficients $a_{*,j}^{(l)},b_{*,j}^{(l)}$ are rational numbers. A standard computation shows that
\begin{align*}
&\int_{0}^{\infty} \frac{t^{l} \; e^{-t}}{(1+e^{t})^{2}} dt=  l! (1-2 (1-2^{-l})\zeta(l+1) + (1-2^{1-l})\zeta(l)),\\
&\int_{0}^{\infty} \frac{t^{l} \; e^{t}}{(1+e^{t})^{2}} dt= l! (1-2^{1-l})\zeta(l).
\end{align*}
From here we see that the values of the zeta function at odd $l=2k+1$ do not contribute to the sum in $(iv)$ if and only if 
\begin{equation}\label{cond}
(2k+1)\cdot(a_{1,j}^{(2k+1)}+a_{2,j}^{(2k+1)})-a_{1,j}^{(2k)}-a_{2,j}^{(2k)}+b_{1,j}^{(2k)}+b_{2,j}^{(2k)}=0.
\end{equation}
The condition (\ref{cond}) is not hard to verify by induction from $j$ to $j+1$. Thus, only the values of the zeta function at even $l$ contribute to the sum.
\hfill $\Box$

\bigskip
\noindent
{\bf References}. 
\begin{itemize}
\item The relationship between the Weil-Petersson volumes and the intersection numbers of $\psi-$classes on $\overline{\mathcal{M}}_{g,n}$  is discussed in \cite{W} and \cite{AC}. 
An explicit formula for the volumes in terms of the intersections of $\psi-$classes was given in \cite{KMZ:w}, cf. also \cite{MZ}.

\item Recursion $({\bf Ia})$ is a special case of Proposition $3.3$ in \cite{LX:higher}. 
Similarly, recursion $({\bf Ib})$ is a simple corollary of Propositions $3.3$ and $3.4$ in 
\cite{LX:higher}.
 
\item For different proofs of $({\bf II})$ see \cite{DN:cone} and \cite{LX:higher}. In terms of the volume polynomial $V_{g,n}(L)$,  recursion
${(\bf II)}$ can be written as follows (\cite{DN:cone}):
$$\frac{\partial V_{g,n+1}}{\partial L_{n+1}}(L_1,\ldots, L_n, 2\pi \sqrt{-1}) = 2\pi \sqrt{-1} (2g-2+ n)V_{g,n}(L_1,\ldots, L_n).$$
When $n=0,$ 
 $$V_{g,1}(2\pi \sqrt{-1})=0,$$
 and
\begin{equation}\label{zero}
\frac{\partial V_{g,1}}{\partial L} (2\pi \sqrt{-1}) = 2\pi \sqrt{-1} (2g-2)V_{g,0}.
\end{equation}

 \item For a proof of $({\bf III})$ see \cite{M:In}; note that $({\bf III})$ applies only when $n>0$ (in case of $n=0,$ formula $(\ref{zero})$ gives the necessary estimates on the growth of $V_{g,0}$). In fact, $({\bf III})$ can be interpreted as a recursive formula for the volumes of  moduli spaces $\mathcal{M}_{g,n}(L)$
 of hyperbolic surfaces of genus $g$ with $n$ geodesic boundary components of lengths $L_{1},\ldots,L_{n}$
that describes a removal of a pair of pants on a surface containing at least one of its boundary
components.  Although $({\bf III})$ is written here in purely
combinatorial terms, it is related to the topology 
of different pant decompositions of a surface, 
cf. also \cite{M:S} and \cite{LX:M}.

 \item When $d_{1}+\ldots+d_{n}=3g-3+n,$ recursion $({\bf III})$ reduces to the Virasoro constraints for the intersection numbers of $\psi_{i}$-classes predicted by Witten  \cite{W}, cf. also \cite{MuS}. 
For different proofs and discussions of these relations see, \cite{Ko:int}, \cite{OP}, \cite{M:JAMS}, \cite{KL:W}, and \cite{EO:WK}.

 \item In this paper, we are mainly interested in the intersection numbers of $\kappa_{1}$ and $\psi_{i}$ classes. For generalizations of $({\bf III})$ to the case of intersection numbers involving higher Mumford's $\kappa$-classes see \cite{LX:higher} , \cite{E:M} and \cite{Ka}.

\end{itemize}

\end{section}

\begin{section}{Asymptotics of intersection numbers when $n$ is fixed}
In this section, we prove Theorem $\ref{11}$.
This theorem implies that there exists $C\in (0,\infty)$ such that 
\begin{equation}\label{first}
\lim_{g\rightarrow \infty} \frac{V_{g,n} \sqrt{g}}{(2g-3+n)!\,(4\pi^2)^{2g-3+n}}=C.
\end{equation}
This result will be generalized in $\ch \ref{error}.$

We recall that for any $n \geq 0$ the results obtained in \cite{M:large} yield 

\begin{equation}\label{f1}
\frac{V_{g,n}}{8\pi^2 g \; V_{g-1,n+1}}= 1+O\left(\frac{1}{g}\right)
\end{equation}
and
\begin{equation}\label{f2}
\frac{V_{g,n}}{ V_{g-1,n+2}}= 1+O\left(\frac{1}{g}\right)
\end{equation}
as $g\to\infty$. The main ingredient of the proof is the following property of the intersection numbers: 

\begin{equation}\label{inequality:basic}
 [\tau_{d_{1}}\ldots\tau_{d_{n}}]_{g,n} \leq [\tau_0^n]_{g,n}=V_{g,n}
 \end{equation}
for any $d=(d_1,\ldots, d_n)$. Moreover, 
\begin{equation}\label{simple}
\frac{[\tau_{d_{1}}\ldots \tau_{d_{n}}]_{g,n}}{V_{g,n}}=1+O\left(\frac{1}{g}\right), 
\end{equation}
as $g \rightarrow \infty.$
\begin{rema}
{\rm The same result holds if $d_{1},\ldots, d_n$ grow slowly with $g$ in such a way that
 $$\frac{d_{1}\ldots d_{n}}{g} \rightarrow 0$$
as $g \rightarrow \infty.$
In particular, $(\ref{simple})$ holds if $d_{i}=O(\log g)$ for each $i=1,\ldots, n$.}
\end{rema}
A stronger statement is formulated below:

\begin{theo}\label{theo:tauk}Let $k,n  \geq 1,$ then
\begin{enumerate}[(i)]
\item
$$ \frac{[\tau_k \,\tau_0^{n-1}]_{g,n}}{V_{g,n}}
=1+\frac{e_{n,k}^{(1)}}{g}+O\left(\frac{1}{g^{2}}\right),$$
as $g \rightarrow \infty$, where
$$e_{n,k}^{(1)}=-\frac{k^2+(n-5/2)k-n/2+3/2}{\pi^2}.$$
\item 
\begin{align*}
  8n [\tau_{k-1}\,\tau_{0}^{n-2}]_{g,n-1}&< [\tau_{k}\,\tau_0^{n-1}]_{g,n} - [\tau_{k+1}\,\tau_0^{n-1}]_{g,n} \leq\\
& \leq  16\; ((n+k) V_{g,n-1}+ k\, V_{g-1,n+1}).
 \end{align*}
 \end{enumerate}
\end{theo}
\begin{rema}\label{general:2:error}
{\rm In general, one can show that for $k\geq 1$
$$\frac{[\tau_{k+1}\,\tau_{j_{1}} \,\ldots
\tau_{j_{s}}\tau_{0}^{n-1-s}]_{g,n}}{[\tau_{k}\,\tau_{j_{1}}\ldots\tau_{j_{s}},
\tau_0^{n-1-s}]_{g,n}}
=1- \frac{2(k+j_1+...j_s)+n-3/2}{\pi^2g}
+O\left(\frac{1}{g^2}\right).$$}
\end{rema}

{\bf  Proof of Theorem $\ref{theo:tauk}$.}
We will need the following simple fact.
Let $\{r_{i}\}_{i=1}^{\infty}$ be a sequence of real numbers and $\{k_{g}\}_{g=1}^{\infty}$ be an increasing sequence of positive integers.
Assume that  for all $g$, and $i$ we have $0 \leq c_{i,g} \leq  c_{i},$ and  
 $\lim_{g\rightarrow \infty} c_{i,g}=c_{i}.$
If $\sum_{i=1}^{\infty} | c_{i} r_{i}| < \infty,$ then 
\begin{equation}\label{fact}
\lim_{g\rightarrow \infty} \sum_{i=1}^{k_{g}} r_{i} c_{i,g} = \sum_{i=1}^{\infty} r_{i} c_{i}.
\end{equation}

To prove part  $(i)$ of the Theorem, it is sufficient to show that
$$\frac{[\tau_{1}\,\tau_{0}^{n-1}]_{g,n}}{[\tau_{0}^n]_{g,n}}=1-  \frac{n}{2\pi^2 g}+O\left(\frac{1}{g^2}\right),$$
and for $k \geq 1$ 
\begin{equation}\label{similar}
\frac{[\tau_{k+1}\,\tau_{0}^{n-1}]_{g,n}}{[\tau_{k}\,\tau_0^{n-1}]_{g,n}}=1- \frac{2k+n-3/2}{\pi^2 g}+O\left(\frac{1}{g^2}\right).
\end{equation}
Here we use the recursive formula $({\bf III})$ to expand the difference $[\tau_{k}\,\tau_{0}^{n-1}]_{g,n}-[\tau_{k+1}\,\tau_{0}^{n-1}]_{g,n}$ in terms of the intersection numbers on $\overline{\mathcal{M}}_{g-1,n+1}$, $\overline{\mathcal{M}}_{g,n-1}$ and $\overline{\mathcal{M}}_{g_{1},n_{1}} \times \overline{\mathcal{M}}_{g_{2},n_{2}}$. For the sake of brevity let us put
\begin{equation}\label{Obs:III}
[\tau_{k}\,\tau_{0}^{n-1}]_{g,n}-[\tau_{k+1}\,\tau_{0}^{n-1}]_{g,n}=\widetilde{{A}}_{k,g,n}+\widetilde{{B}}_{k,g,n}+\widetilde{{C}}_{k,g,n},
\end{equation}
where $\widetilde{{A}}_{k,g,n}$, $\widetilde{{B}}_{k,g,n}$, and $\widetilde{{C}}_{k,g,n}$ are the terms corresponding to $(\ref{A}),$ $(\ref{B})$ and $(\ref{C})$ respectively.
We will evaluate these terms separately.

\medskip
\noindent
{\bf 1. }{\em Contribution from} ($\ref{A}$).
By $({\bf III})$ , the numbers $$[\tau_{k-1}\,\tau_{0}^{n-2}]_{g,n-1}, \ldots, [\tau_{3g+n-4}\,\tau_{0}^{n-2}]_{g,n-1}$$ contribute to $[\tau_{k}\,\tau_{0}^{n-1}]_{g,n}$ and $[\tau_{k+1}\,\tau_{0}^{n-1}]_{g,n}$.
In fact, it is easy to check that 
$$\widetilde{{A}}_{k,g,n}= 8 \,a_{0}\; [\tau_{k-1}\,\tau_{0}^{n-2}]_{g,n-1}\; +\;8\sum_{i=1}^{3g-3+n-k} (a_{i+1}-a_{i}) [\tau_{k-1+i}\,\tau_{0}^{n-2}]_{g,n-1}.$$
The term $[\tau_{k-1}\,\tau_{0}^{n-2}]_{g,n-1}$ is non-zero only when $k\geq 1.$
In order to calculate the asymptotic behavior of $\widetilde{{A}}_{k,g,n}/ V_{g,n-1}$ we simply apply $(\ref{simple})$ and $(\ref{fact})$.
In view of Lemma $\ref{aa}$, we get that 
$$\lim_{g \rightarrow \infty} \frac{\widetilde{{A}}_{k,g,n}}{V_{g,n-1}}= 8 (n-1)\;\left (a_{0} \delta+ \sum_{i=0}^{\infty} (a_{i+1}-a_{i})\right)= 8(n-1) (1/2 \; \delta+1/2),$$
where $\delta=0$ when $k=0,$ and otherwise $\delta=1.$ Thus
\begin{equation}\label{yek}
 \frac{\widetilde{{A}}_{k,g,n}}{V_{g,n-1}}= 8(n-1) (1/2 \; \delta+1/2)+ O\left(\frac{1}{g}\right)
\end{equation}
as $g\to\infty$. On the other hand, Lemma \ref{aa} also implies that
\begin{equation}\label{1111}
\widetilde{{A}}_{k,g,n} \leq  8(n-1) V_{g,n-1}.
\end{equation}

\medskip
\noindent
 {\bf 2. }{\em Contribution from} ($\ref{B}$).
Similarly, by $({\bf III})$, the numbers $[\tau_{i}\,\tau_{j}\, \tau_{0}^{n-1}]_{g-1,n+1}$ contribute to $[\tau_{k}\,\tau_{0}^{n-1}]_{g,n}$ (resp. to $[\tau_{k+1}\,\tau_{0}^{n-1}]_{g,n}$) whenever  $i+j\geq k-2$ (resp. $i+j \geq k-1$). To simplify the notation, let 
$$T_{m,g,n}= \sum_{i+j=m} [\tau_{i}\,\tau_{j}\,\tau_{0}^{n-1}]_{g-1,n+1}$$
(we assume  $T_{m,g,n}=0$ for $m <0$).
Then
$$\widetilde{{B}}_{k,g,n}=16\,( a_{0} T_{k-2,g,n}+ (a_{1}-a_{0}) T_{k-1,g,n}+\ldots+ (a_{i+1}-a_{i}) T_{k-1+i,g,n}+\ldots).$$
Note that by $(\ref{simple})$ as $g \rightarrow \infty,$
$$ \lim_{g\rightarrow \infty} \frac{T_{m,g,n}}{V_{g-1,n+1}}=m+1.$$
Now since $a_{0}=1/2,$ Lemma $\ref{aa}\; (iii)$, together with $(\ref{simple})$ and $(\ref{fact})$, implies for $k>1$

\begin{align}
\lim_{g\rightarrow \infty} \frac{\widetilde{{B}}_{k,g,n}}{V_{g-1,n+1}}&=16\left( a_{0}(k-1)+ \sum_{i=0}^{\infty} (i+k) ( a_{i+1}-a_{i})\right)= \nonumber\\
&=16\left(\frac{k-1}{2}+ \frac{k}{2}+1/4\right)=16( k-1+1/4),\label{do}\\
\lim_{g\rightarrow \infty} \frac{\widetilde{{B}}_{1,g,n}}{V_{g-1,n+1}}&=16\sum_{i=0}^{\infty} (i+1) (a_{i+1}-a_{i})=4,\nonumber
\end{align}
and we also have
\begin{equation}\label{222}
\widetilde{{B}}_{k,g,n} \leq  16 k V_{g-1,n+1}.
\end{equation}

\medskip
\noindent
{\bf 3. }{\em Contribution from} ($\ref{C}$).
By the results obtained in \cite{M:large}
$$ \sum_{\sumind{I \amalg J=\{2,\ldots,n\}}{0\leq g' \leq g}} \frac{ V_{g',|I|+1} \cdot V_{g-g',|J|+1}}{V_{g,n}}= O\left(\frac{1}{g^2}\right).$$
Put 
$$S_{k_1,k_2,g,n}= \sum_{\sumind {I \amalg J=\{2,\ldots,n\}}{0\leq g' \leq g}} \; \left[\tau_{k_{1}} \prod_{i\in I } \tau_{i}\right]_{g',|I|+1} \cdot \quad\left[ \tau_{k_{2}} \prod_{i\in J} \tau_{i}\right]_{g-g',|J|+1}.$$
Note that by $(\ref{inequality:basic})$ and recursion $({\bf Ia})$,
$$ S_{k_{1},k_{2},g,n} \leq \sum_{\sumind{I \amalg J=\{2,\ldots,n\}}{0\leq g' \leq g}} V_{g',|I|+1} \cdot V_{g-g', |J|+1} \leq V_{g,n-1}$$
Therefore, the contribution from the term $(\ref{C})$ in $({\bf III})$ satisfies
\begin{align*}
0& \leq \widetilde{{C}}_{k,g,n} \leq 16\;\sum_{i=0}^\infty (a_{i}-a_{i+1}) \sum_{k_{1}+k_{2}=i+k} S_{k_{1},k_{2},g,n} \leq\\
& \leq 16\sum_{i=0}^\infty (i+k) (a_{i}-a_{i+1})\; V_{g,n-1} \leq 16\,(1/4+k/2)\;V_{g,n-1}.
 \end{align*}
Using $(\ref{fact})$, as in the cases  {\bf 1} and  {\bf 2} considered above, we see that the contribution from the term $(\ref{C})$ in $({\bf III})$ becomes small as $g\to\infty$:
\begin{equation}\label{se}
 \frac{\widetilde{{C}}_{k,g,n}}{V_{g,n}}=O\left(\frac{1}{g^{2}}\right),
 \end{equation}
and 
\begin{equation}\label{333}
\widetilde{{C}}_{k,g,n} \leq 16 (1+k) \,V_{g,n-1}.
\end{equation}

\medskip
Now, in view of $(\ref{f1})$, $(\ref{f2})$ and $(\ref{simple})$, equations $(\ref{yek})$, $(\ref{do})$ and $(\ref{se}))$ imply that for $k\geq 1$
$$1-\frac{[\tau_{k+1}\,\tau_{0}^{n-1}]_{g,n}}{[\tau_{k}\,\tau_{0}^{n-1}]_{g,n}}= \frac{2k+n-3/2}{\pi^2}\cdot\frac{1}{g}+O\left(\frac{1}{g^2}\right)$$
and 
$$1-\frac{[\tau_1\,\tau_0^{n-1}]_{g,n}}{[\tau_0^n]_{g,n}}= \frac{n}{2\pi^2}\cdot\frac{1}{g}+O\left(\frac{1}{g^2}\right).$$
Finally, the inequalities $(\ref{1111}),$ $(\ref{222})$ and $(\ref{333})$ imply part $(ii)$ of the Theorem.
\hfill $\Box$

\noindent
{\bf Proof of Theorem \ref{11}}.
We start with proving $(\ref{110})$. From $({\bf II})$, 
$$ \frac{2(2g-2+n) V_{g,n}}{V_{g,n+1}}= 
 \sum_{l=1}^{3g-2+n} \frac{(-1)^{l-1}l\, \pi^{2l-2}}{(2l+1)!}\cdot \frac{[ \tau_{l},\tau_{0}^n]_{g,n+1}}{V_{g,n+1}}.$$
Differentiating $t^{-1}\sin t$ and putting $t=\pi$ we get
$$ \sum_{l=1}^{\infty} \frac{(-1)^{l-1} l \,\pi^{2l-2}}{(2l+1)!}= \frac{1}{2\pi^2}.$$
Now we can use $(\ref{fact})$, and Theorem $\ref{theo:tauk}$ in order to calculate the error term in $(2g-2+n) V_{g,n}/V_{g,n+1}-1/4\pi^{2}.$ 
Clearly,
\begin{equation}\label{Obs:II}
 \frac{(2g-2+n) V_{g,n}}{V_{g,n+1}}-\frac{1}{4\pi^2} = 
 \sum_{l=1}^{3g-2+n} \frac{(-1)^{l-1} l\, \pi^{2l-2}}{(2l+1)!} \cdot\left(\frac{[ \tau_{l}\,\tau_{0}^n]_{g,n+1}}{V_{g,n+1}}-1\right).
\end{equation}
Then Theorem \ref{theo:tauk}, $(i)$ implies that 
\begin{align*}
\frac{(2g-2+n) V_{g,n}}{V_{g,n+1}}&-\frac{1}{4\pi^{2}}=\\
= &- \sum_{l=1}^{\infty}\frac{ (-1)^{l-1} (l^2+(n-3/2)l-n+1)\,l\, \pi^{2l-2}}{2 (2l+1)! \pi^2}\cdot\frac{1}{g}+O\left(\frac{1}{g^{2}}\right).
 \end{align*}
On the other hand, 
$$ \sum_{l=1}^{\infty} \frac{(-1)^{l-1} \,l\, (l^2+(n-3/2)l-n/2+1) \pi^{2l}}{(2l+1)!}= -\frac{4n+\pi^2-8}{8\pi^2}.$$ 
Hence, 

\begin{equation}\label{eq:n} 
\frac{4\pi^2 (2g-2+n)\; V_{g,n}}{V_{g,n+1}}= 1+\frac{4n+\pi^2-8}{4\pi^2}\cdot\frac{1}{g}+O\left(\frac{1}{g^{2}}\right),
\end{equation}
and
$$  \frac{8\pi^2g\; V_{g,n}}{V_{g,n+1}}=1+\left(\left(\frac{1}{\pi^2}-\frac{1}{2}\right)n+\frac{5}{4}-\frac{2}{\pi^2}\right)\cdot\frac{1}{g}+ O\left(\frac{1}{g^{2}}\right).$$

We proceed with proving $(\ref{111})$.
First, we will check this estimate when $n \geq 2$, that is,
$$ \frac{V_{g-1,n+4}} {V_{g,n+2}}=1-\frac{2n+1}{ \pi^2}\cdot\frac{1}{ g}+O\left(\frac{1}{g^2}\right).$$  
From the recursion $({\bf Ia})$ we get
\begin{align}\frac{V_{g-1,n+4}}{V_{g,n+2}}&=
\frac{[ \tau_{1}\,\tau_{0}^{n+1}]_{g,n+2} }{V_{g,n+2}}-\nonumber\\ \label{Obs:I}
&-\frac{6}{V_{g,n+2}} \sum_{\sumind{g_1+g_2=g}{I\amalg J=\{1,\ldots, n\}}} V_{g_{1},|I|+2} \cdot V_{g_{2},|J|+2}.
\end{align}
On the other hand, by Theorem $\ref{theo:tauk}$ we have for $k=1$
$$\frac{[\tau_{1}\,\tau_{0}^{n+1}]_{g,n+2}}{[\tau_{0}^{n+2}]_{g,n+2}}=1-  \frac{n+2}{2\pi^2}\cdot\frac{1}{g}+O\left(\frac{1}{g^2}\right),$$
and $(\ref{f2})$ implies
\begin{align*}
\frac{1}{V_{g,n+2}} \sum_{\sumind{g_{1}+g_{2}=g}{ I\amalg J=\{1,\ldots, n\}} } V_{g_{1},|I|+2} \cdot V_{g_{2},|J|+2}
&= 2 n\; \frac{V_{g,n+1}}{V_{g,n+2}} +O\left(\frac{1}{g^{2}}\right)\\
&= \frac{2 n}{8\pi^{2}} \cdot\frac{1}{g}+O\left(\frac{1}{g^{2}}\right).
\end{align*}
Hence
$$ \frac{V_{g-1,n+4}}{V_{g,n+2}}=1- \left(\frac{n+2}{2\pi^2}+ \frac{12\; n}{8\pi^{2}}\right)\cdot\frac{1}{g}+ O\left(\frac{1}{g^{2}}\right)= 1-\frac{2n+1} {\pi^2 g}+O\left(\frac{1}{g^{2}}\right). $$
The remaining cases $n=0,\,1$ follow from (\ref{110}) and (\ref{111}) for $n\geq2$. For instance, if $n=1$
$$
\frac{V_{g-1,3}}{V_{g,1}}=\frac{V_{g-1,4}}{V_{g,2}}\cdot\frac{V_{g,2}}{V_{g,1}}\cdot\frac{V_{g-1,3}}{V_{g-1,4}}
=1+\frac{1}{\pi^2}\cdot\frac{1}{g}+O\left(\frac{1}{g^2}\right).
$$
The case $n=0$ can be treated similarly.
\hfill $\Box$

Theorem \ref{11} immediately implies Corollary $\ref{c}$ about the asymptotic behavior of the ratio
$V_{g+1,n}/V_{g,n}$ (see Introduction). An important consequence of Corollary $\ref{c}$, explained in Remark $\ref{rem2}$, is formula $(\ref{first}$) announced at the beginning of this Section.
As a byproduct of this statement we also get the following estimate that we will need later:
\begin{lemm}\label{estimate:upperbound}
Fix $n_{1}, n_{2},s \geq 0.$ Then
\begin{equation}
\sum_{\sumind{g_{1}+g_{2}=g}{2g_{i}+n_{i} \geq s,\; i=1,2}} V_{g_1,n_1} \cdot V_{g_{2},n_{2}}= O\left(\frac{V_{g,n_{1}+n_{2}}}{g^{s}}\right).
\end{equation}
\end{lemm}
\end{section}

\begin{section}{Error terms in the asymptotics expansions}\label{error}
In this section, we prove Theorem $\ref{theo:main:fixn}$ using the following results:

\begin{theo}\label{theo:general:1} We have the following asymptotic expansions as $g \rightarrow \infty:$
\begin{enumerate}[(i)]
\item
Given the integers $n, s\geq 1,$ and ${d}=(d_{1},\ldots, d_{n})$, there exist $e_{n,{d}}^{(1)},\ldots e_{n,{d}}^{(s-1)}$ independent of $g$ such that 
\begin{equation}
 \frac{[\tau_{d_1} \ldots \tau_{d_n}]_{g,n}}{V_{g,n}}= 1+\frac{e_{n,{d}}^{(1)}}{g}+\ldots+\frac{e_{n,{d}}^{(s-1)}}{g^{s-1}}
 +O\left(\frac{1}{g^{s}}\right).\label{e}
 \end{equation}
\item Given $n\geq 0,$ $s\geq 1$, there exist $a_{n}^{(i)},\;b_{n}^{(i)}\; i=1,\ldots , s-1,$ independent of $g$ such that
\begin{align}
\frac{4\pi^2 (2g-2+n)V_{g,n}}{V_{g,n+1}} = &1+ \frac{a_{n}^{(1)}}{g}+\ldots+\frac{a_{n}^{(s-1)}}{g^{s-1}}
+O\left(\frac{1}{g^{s}}\right),\label{theo:n}\\
\frac{V_{g,n}} {V_{g-1,n+2}}=&1+ \frac{b_{n}^{(1)}}{g}+\ldots+\frac{b_{n}^{(s-1)}}{g^{s-1}}+O\left(\frac{1}{g^{s}}\right).\label{theo:g}
\end{align}
\end{enumerate}
\end{theo}

The coefficients of the above asymptotic expansions (\ref{e})--(\ref{theo:g}) can be characterized more precisely:
\begin{theo}\label{theo:general:2}
We have
\begin{enumerate}[(i)]
\item For any fixed $n$ and $d$ the coefficient $e_{n,{d}}^{(i)}$ is a  polynomial in ${\mathbb Q}[\pi^{-2}, \pi^{2}]$ of degree at most $i.$
\item Each $a_{n}^{(i)}$ and $b_{n}^{(i)}$ is a polynomial in  ${\mathbb Q}[\,n,\,\pi^{-2}, \pi^2]$ of degree $i$ in $n$ and of degree at most $i$ in $\pi^{-2}$ and $\pi^2$.
\end{enumerate}
\end{theo}

\begin{rema}
{\rm In the simplest case $[\tau_0\tau_k]_{g,2}/V_{g,2}$ we have the following expansions:
\begin{align*}
\frac{[\tau_0\tau_1]_{g,2}}{V_{g,2}}
&=1-\frac{1}{\pi^2g}+\left(\frac{1}{64}-\frac{5}{6\pi^2}+\frac{1}{\pi^4}\right)\cdot\frac{1}{g^2}+O\left(\frac{1}{g^3}\right),\\
\frac{[\tau_0\tau_2]_{g,2}}{V_{g,2}}
&=1-\frac{7}{2\pi^2g}+\left(\frac{1}{64}-\frac{13}{6\pi^2}+\frac{1}{\pi^4}\right)\cdot\frac{1}{g^2}+O\left(\frac{1}{g^3}\right),\\
\frac{[\tau_0\tau_k]_{g,2}}{V_{g,2}}
&=1-\frac{2k^2-k+1}{2\pi^2g}+\left(\frac{k^4}{2\pi^4}-\frac{13k^3}{6\pi^4}-\left(\frac{1}{2\pi^2}-\frac{27}{8\pi^4}\right)\cdot k^2\right.\\
&+\left.\left(\frac{1}{24\pi^2}-\frac{59}{24\pi^4}\right)\cdot k+\left(\frac{1}{64}-\frac{1}{4\pi^2}+\frac{19}{8\pi^4}\right)\right)\cdot\frac{1}{g^2}+O\left(\frac{1}{g^3}\right).
\end{align*}
We see that no (positive) powers of $\pi^2$ appear in these expansions.
The term of order $1/g^2$ is a polynomial in $k$ of degree 4 for $k\geq 3$ (computed numerically). However,  for $k=1,2$ the general formula is off by $\frac{1}{8\pi^2}+\frac{5}{8\pi^4}$
and $\frac{5}{8\pi^4}$ for $k=1$ and $k=2$ respectively. This is a manifestation of the ``boundary effect" in recursion ({\bf III}). These results can be proved using 
Remark $\ref{general:2:error}$.}
\end{rema}

\begin{rema}\label{rem33}
{\rm Note that a result similar to $(\ref{theo:n})$ and $(\ref{theo:g})$ holds for the inverse ratios $V_{g,n+1}/ (8\pi^2 g V_{g,n})$ and $V_{g-1,n+2}/V_{g,n}.$
This is because of the following simple fact. Let $\{w_{g}\}_{g=1}^{\infty}$ be a sequence of the form
$$w_{g}=1+ \frac{u_1}{g}+\ldots+\frac{u_{s-1}}{g^{s-1}}+O\left(\frac{1}{g^{s}}\right),$$
then
$$ \frac{1}{w_{g}}= 1+ \frac{v_{1}}{g}+\ldots+\frac{v_{s-1}}{g^{s-1}}+O\left(\frac{1}{g^{s}}\right),$$
where each $v_{i}$ is a polynomial in $u_{1},\ldots, u_{i}$ with integer coefficients. Moreover, if $u_{i}$ is a polynomial of degree $m_{i}$ in $n$, then
$v_k$ is a polynomial of degree at most $\max_{i+j=k}(m_{i}+m_{j})$ in $n$.}
\end{rema}

Let us first outline some general ideas underlying the proofs of Theorems $\ref{theo:general:1}$ and $\ref{theo:general:2}$.
All proofs are by induction in $s$ and are similar to each other. We basically follow the same steps as in the course of proving Theorems $\ref{theo:tauk}$ and $\ref{11}$. 

\begin{rema}\label{rem3}

{\rm Let $f: {\mathbb Z}_+\rightarrow {\mathbb R}$ be a function such that $\lim_{g \rightarrow \infty} f(g)$ exists. We say that $f$ has an expansion up to $O(1/g^{s})$ if
there exist $e_{0},\ldots,e_{s-1} \in {\mathbb R}$ so that 
$$f(g)=e_{0}+\frac{e_{1}}{g}+\ldots+ \frac{e_{s-1}}{g^{s-1}}+O\left(\frac{1}{g^{s}}\right).$$
Note that if $f_{1},\ldots,f_{k}$ all have expansions of order $s,$ the expansion of the product $f_{1}\cdots f_{k}$ up to $1/g^{s}$ can easily be calculated in terms of the expansions of  $f_{i}.$ 

In this section, we are interested in the expansions of the ratios $\frac{V_{g,n}}{V_{g-1,n+2}}$, 
$\frac{4\pi^2(2g-2+n) V_{g,n}}{V_{g,n+1}}$
and 
$\frac{[\tau_{d_1}\ldots\tau_{d_{n}}]_{g,n}}{V_{g,n}}. $ Some remarks are in order:

\medskip
\noindent {\bf 1.}
In general, given $g', n'$, in order to obtain the expansion of $\frac{V_{g-g',n-n'}}{V_{g,n}}$ up to $O(1/g^s),$ it is enough to know the expansions of
$\frac{V_{g,k}}{V_{g-1,k+2}}$ and 
$\frac{4\pi^2 (2g-2+k) V_{g,k}}{V_{g,k+1}}$ up to $O(1/g^{s-2g'-n'});$ this is simply because $$\frac{V_{g-g',n-n'}}{V_{g,n}}= \prod_{j=-n'+1}^{2g'} \frac{4 \pi^2 (2g-2g'+n-j+1)V_{g-g',n+j-1}}{V_{g-g',n+j}} \cdot$$ 
 \begin{align}\label{subsurface} 
\cdot \prod_{j=1}^{g'} \frac{V_{g-j,n+2j}}{V_{g-j+1,n+2j-2}} \cdot  \prod_{j=-n'+1}^{2g'} \frac{1}{(4 \pi^2 (2g-2g'+n-j+1))}.
  \end{align}

\medskip
\noindent 
{\bf 2.} Following $({\bf Ia})$ and $({\bf II})$, the expansion of $\frac{4\pi^2(2g-2+n) V_{g,n}}{V_{g,n+1}}$  up to $O(1/g^s)$ can be 
written explicitly in terms of the expansion of $\frac{[\tau_{l}\tau_{0}^{n}]_{g,n+1}}{V_{g,n+1}}$ up to $O(1/g^{s})$; see (\ref{Obs:II}). 
Similarly, by $(\ref{Obs:III})$ and $(\ref{subsurface})$ the expansion of $\frac{V_{g,n}}{V_{g-1,n+2}}$ up to $O(1/g^{s})$
can be written in terms of the expansions of $\frac{V_{g,n_1}}{V_{g-1,n_1+2}}$ and $\frac{4\pi^2(2g-2+n_1) V_{g,n_1}}{V_{g,n_1+1}}$ up to 
$O(1/g^{s-1}),$ and the expansion of  $\frac{[\tau_{1}\tau_0^{n-1}]_{g,n}}{V_{g,n}}$ up to $O(1/g^s).$

\medskip
\noindent 
{\bf 3.}  In view of $({\bf III}),$ the expansion of $\frac{[\tau_{d_1}\ldots\tau_{d_{n}}]_{g,n}}{V_{g,n}}$ up to $O(1/g^{s})$ 
can be written in terms of the expansions of $\frac{[\tau_{c_{1}}\ldots\tau_{c_{m}}]_{g,m}}{V_{g,m}}$
up to $O(1/g^{s-1})$ and the expansions of  $\frac{V_{g,n}}{V_{g-1,n+2}}$ and $\frac{4\pi^2(2g-2+n) V_{g,n}}{V_{g,n+1}}$ 
up to $O(1/g^{s-1}).$ Actually, for our purposes
 it will be enough to obtain the expansion of 
 $\frac{[\tau_{d_{1}+1}\,\ldots \tau_{d_n}]_{g,n}}{[\tau_{d_{1}}\,\ldots \tau_{d_n}]_{g,n}}$ up to 
$O(1/g^s).$ As in $(\ref{Obs:III})$, we put 
$$[\tau_{d_{1}}\,\ldots \tau_{d_n}]_{g,n}-[\tau_{d_{1}+1}\,\ldots \tau_{d_n}]_{g,n}=\widetilde{{A}}_{{d},g,n}+\widetilde{{B}}_{{d},g,n}+\widetilde{{C}}_{{d},g,n},$$
where $\widetilde{A}_{d,g,n}$, $\widetilde{{B}}_{{d},g,n}$, and $\widetilde{{C}}_{{d},g,n}$ are the terms corresponding to $(\ref{A}),$ $(\ref{B})$ and $(\ref{C})$. Put
\begin{equation}\label{eq:main}
\frac{[\tau_{d_{1}}\,\ldots \tau_{d_n}]_{g,n}-[\tau_{d_{1}+1}\,\ldots \tau_{d_n}]_{g,n}}{V_{g,n}}= S_{1}+S_{2}+S_{3}\;,
\end{equation}
where
\begin{align*}
S_1=&\frac{1}{4 \pi^2 (2g-3+n)}\cdot \frac{4\pi^2 (2g-3+n)V_{g,n-1}}{V_{g,n}}\cdot\frac{\widetilde{{A}}_{{d},g,n}}{V_{g,n-1}},\\
S_2= &\frac{1}{4\pi^2 (2g-3+n)}\cdot \frac{4 \pi^2 (2g-3+n)V_{g-1,n+1}}{V_{g-1,n+2}}\cdot \frac{V_{g-1,n+2}}{V_{g,n}}
\cdot\frac{\widetilde{{B}}_{{d},g,n}}{V_{g-1,n+1}},\\
S_3=&\frac{\widetilde{{C}}_{{d},g,n}}{V_{g,n}}.
\end{align*}
Similar to the case ${\bf 1}$ in the proof of Theorem $\ref{theo:tauk}$, we have
\begin{equation}\label{eq:aa}
\frac{\widetilde{{A}}_{{d},g,n}}{V_{g,n-1}}=8\;\sum_{i=0}^{3g-3+n-d_{1}} (a_{i+1}-a_{i}) \frac{[\tau_{d_{1}+d_{j}+i-1}\tau_{d_{2}}\ldots\widehat{\tau_{d_{j}}}\ldots \tau_{d_{n}}]_{g,n-1}}{V_{g,n-1}} 
\end{equation}
(the hat means that the corresponding entry is omitted, and $a_{-1}=0$).
The case {\bf 2} of the same proof now reads
\begin{equation}\label{eq:bb}
\frac{\widetilde{{B}}_{{d},g,n}}{V_{g-1,n+1}}=16\,(a_{0} T_{k-2,g,n}+ (a_{1}-a_{0}) T_{k-1,g,n}+\ldots+(a_{i+1}-a_{i}) T_{k-1+i,g,n}+\ldots),
\end{equation}
where 
$$T_{m,g,n}=\frac{\sum_{i+j=m} [\tau_{i}\tau_{j}\tau_{d_{2}}\ldots \tau_{d_{n}}]_{g-1,n+1}}{V_{g-1,n+1}}.$$
Similarly, according to $(\ref{C})$, each term in $\widetilde{{C}}_{{d},g,n}$ has the form
\begin{equation}\label{eq:cc}
\sum_{\genfrac{}{}{0pt}{}{k_{1}+k_{2}=}{=i+d_{1}-2}} (a_{i+1}-a_{i}) \; \left[\tau_{k_{1}} \prod_{i\in I } \tau_{d_{i}}\right]_{g',|I|+1} \cdot \left[\tau_{k_{2}} \prod_{i\in J} \tau_{d_{i}}\right]_{g-g',|J|+1},
\end{equation}
where $I \amalg J=\{2,\ldots,n\},$ and $0\leq g' \leq g.$
In order to obtain the expansions of $S_{1}$ and $S_2$, we can use the expansions of ratios 
$\frac{[\tau_{c_{1}}\ldots\tau_{c_{n-1}}]_{g,n-1}}{V_{g,n-1}}$ and 
$\frac{[\tau_{c_{1}}\ldots\tau_{c_{n+1}}]_{g-1,n+1}}{V_{g-1,n+1}}$ up to $O(1/g^{s-1})$.  
What concerns the term $S_3$, by  Lemma $\ref{estimate:upperbound}$ each product
$$\left[\tau_{k_{1}} \prod_{i\in I } \tau_{d_{i}}\right]_{g',|I|+1} \cdot \left[\tau_{k_{2}} \prod_{i\in J} \tau_{d_{i}}\right]_{g-g',|J|+1}$$
is of order $O(1/g^{s+1})$ unless either $2g'-1+|I| <s $ or $2g-2g'+|J|-1 < s$. In these cases we apply 
$(\ref{subsurface})$ to obtain the expansion of $S_{3}$ up to $O(1/g^s).$
Then we can use the expansions of $\frac{4\pi^2(2g-2+n)V_{g-1,n+2}}{V_{g-1,n+3}}$ (for $S_1$) and 
$\frac{V_{g,n}}{V_{g-1,n+2}},\;\frac{4\pi^2(2g-3+n) V_{g-1,n+1}}{V_{g-1,n+2}}$ (for $S_2$), all up to $O(1/g^{s-1}),$ to get the expansion of  $\frac{[\tau_{d_{1}}\,\ldots \tau_{d_n}]_{g,n}-[\tau_{d_{1}+1}\,\ldots \tau_{d_n}]_{g,n}}{V_{g,n}}$ up to $O(1/g^s),$
which will complete the inductive step.}
\end{rema}

\begin{rema}\label{rneed}
{\rm We will need the following basic facts to prove Theorems $\ref{theo:general:1}$ and $\ref{theo:general:2}$:
\begin{enumerate}
\item For any $k\geq 0$ the sum
$$\sum_{l=0}^{\infty} \frac{(-1)^{l} l^{k} \pi^{2l}}{(2l+1)!}$$ 
is a polynomial in $\pi^{2}$ of degree at most $2[k/2]$ with rational coefficients 
(this can be easily seen by expanding $\frac{\sin x}{x}$ in the Taylor series, differentiating it and putting $x=\pi$). 
\item For a polynomial $p(x)=\sum_{j=1}^{m} b_{j} x^{j}$ of degree $m$,
$$\tilde{p}(x)= \sum_{i=1}^{\infty} (a_{i+1}-a_{i}) p(x+i)$$
is again a polynomial of degree $m$. The coefficient of $\tilde{p}(x)$ at $x^{j}$ 
is equal to $\sum_{j+r \leq m} \genfrac{(}{)}{0pt}{}{j+r}{j} b_{j+r} \cdot A(r)$ 
where $A(r)=\sum_{i=1}^{\infty} i^{r} (a_{i+1}-a_{i}).$
\item Since $\psi_i$ and $\kappa_1=\frac{[\omega_{g,n}]}{2\pi^2}$ are rational classes (i.e., belong to $H^{2}(\overline{\mathcal{M}}_{g,n},{\mathbb Q}),$ cf. \cite{W:H}, \cite{AC}), $$[\tau_{d_{1}}\ldots \tau_{d_{n}}]_{g,n} \in {\mathbb Q}\cdot \pi^{6g-6+2n-2|d|},$$
where $|d|=\sum_{i=1}^{n} d_{i}$. In particular, $V_{g,n}=[\tau_0\ldots\tau_0]_{g,n}$ is a rational multiple of $\pi^{6g-6+2n}$ (\cite{W:H}, see also \cite{M:In} for a different point of view).
\item The function $S_{m}(n)=\sum_{i=1}^{n} i^{m}$ is a polynomial in $n$ of degree $m+1$ with rational coefficients (Faulhaber's formula).
\end{enumerate}}
\end{rema}

\noindent
{\bf Proof of Theorem $\ref{theo:general:1}$}. First, we use $({\bf Ia}),$ $({\bf II})$ and $({\bf III})$ to prove the existence of $e_{n,d}^{(s)}$, $a_{n}^{(s)}$ and $b_{n}^{(s)}$. This is similar to what we did in the proofs of $(\ref{110})$ and $(\ref{111}).$
In fact, instead of $(i)$ we will prove a stronger statement: namely, there exist polynomials $Q_{n}^{{(s)}}(d_1,\ldots,d_n)$ 
and $q_{n}^{{(s)}}(d_1,\ldots,d_n)$ in variables $d_1,\ldots,d_n$ of degrees $s+1$ and $s$ respectively such that for any $d=(d_{1},\ldots,d_{n})$
\begin{equation}\label{theo:ineq1}
\left|\frac{[\tau_{d_1} \ldots \tau_{d_n}]_{g,n}}{V_{g,n}}-1-\frac{e_{n,{d}}^{(1)}}{g}-\ldots-\frac{e_{n,{d}}^{(s)}}{g^{s}}\right| 
\leq \frac{Q_{n}^{(s)}(d_{1},\ldots,d_{n})}{g^{s+1}},
\end{equation}
and
\begin{equation}\label{theo:ineq2}
|e_{n,{d}}^{(s)}| \leq q_{n}^{s}(d_{1},\ldots,d_{n}).
\end{equation}
These formulas follow from the two claims below:

\medskip 
\noindent
{\bf Claim $1$}:  {\em Formulas $(\ref{theo:ineq1})$ and $(\ref{theo:ineq2})$ for $s=r$ and formulas $(\ref{theo:n}),$ $(\ref{theo:g})$ 
for $s<r$ imply $(\ref{theo:n})$ and $(\ref{theo:g})$ for $s=r$.}

\medskip
In fact,  from
$(\ref{Obs:II})$, $(\ref{theo:ineq1})$ and $(\ref{theo:ineq2})$ for $d=(l,0,\ldots,0)$ we have
\begin{equation}\label{f11}
\frac{(2g-2+n)V_{g,n}} {V_{g,n+1}}-\frac{1}{4\pi^2}= \frac{a_{n}^{(1)}}{g}+\ldots+\frac{a_{n}^{(s)}}{g^{s}}+O\left(\frac{1}{g^{s+1}}\right),
\end{equation}
where
\begin{equation}\label{a:n:s}
a_{n}^{(s)}=  \sum_{l=1}^{\infty} \frac{(-1)^{l-1} l\, \pi^{2l-2}}{(2l+1)!} e_{n,l}^{(s)},
\end{equation}
with $d=(l,0,\dots,0).$ The existence of $a_{n}^{(s)}$ is guaranteed by the estimate
$$ \sum_{l=N}^{\infty} \frac{(-1)^{l-1} l^{k+1}\, \pi^{2l-2}}{(2l+1)!}=O(e^{-N})$$ 
valid for any $k\geq 0$.

Similarly, we can use $(\ref{Obs:I})$ to evaluate the error term in the expansion of  $\frac{V_{g-1,n+4}}{V_{g,n+2}}-1$. 
In this case we apply $(\ref{theo:ineq1})$ with $d=(1,0,\ldots,0)$. Note that by Lemma $\ref{estimate:upperbound}$
\begin{align}
\frac{6}{V_{g,n+2}} &\sum_{\sumind{g_1+g_2=g}{I\amalg J=\{1,\ldots, n\}}} V_{g_{1},|I|+2} \cdot V_{g_{2},|J|+2}= \nonumber\\
=&\sum_{2j+i+2 \leq s} \genfrac{(}{)}{0pt}{}{n}{i} \frac{V_{g-j,n+2-i}}{V_{g,n+2}} \times V_{j,i+2}+
O(\frac{1}{g^{s+1}}).\label{next}
\end{align}
We can now use $(\ref{subsurface})$, and together with 
$(\ref{theo:n})$ and $(\ref{theo:g})$ for $s=r$ this yields the expansion of $\frac{V_{g-j,n+2-i}}{V_{g,n+2}}$ up to $O(1/g^{r+1}).$ 

\medskip
\noindent
{\bf Claim $2$}. {\em Fromulas $(\ref{theo:n})$, $(\ref{theo:g}),$ $(\ref{theo:ineq1})$ and $(\ref{theo:ineq2})$ for $s<r$ imply $(\ref{theo:ineq1})$ and $(\ref{theo:ineq2})$ for $s=r$. }

\medskip
According to $(\ref{eq:main}),$ we need to evaluate the contributions from the term $S_{1},$ $S_{2}$ and $S_{3}$ up to $O(1/g^r)$. In view of $(\ref{eq:aa}),$ $(\ref{eq:bb})$ and $(\ref{eq:cc})$ we can use $(\ref{theo:ineq1})$ for $s=r-1$ to obtain the expansions of  $\frac{\widetilde{{A}}_{{d},g,n}}{V_{g,n-1}}$ and$\frac{\widetilde{{B}}_{{d},g,n}}{V_{g-1,n+1}}$ up to $O(1/g^{r-1}).$ Formula $(\ref{eq:aa})$ now takes the form
\begin{equation}\label{tf}
\frac{\widetilde{{A}}_{{d},g,n}}{V_{g,n-1}}=8\,\sum_{i=0}^{\infty} (a_{i+1}-a_{i}) 
\left(1+\frac{e_{n,{d(i)}}^{(1)}}{g}+\ldots+\frac{e_{n,{d'}}^{(r-1)}}{g^{r-1}}+E_{d',r}\right),
\end{equation}
where $d(i)=d_{1}+d_{j}+i-1, d_{2},\ldots\widehat{d_{j}}\ldots d_{n},$ and $E_{d',r} \leq \frac{Q_{n}^{r-1}(d')}{g^{r}}$. Note that by Lemma $\ref{aa}$, $(ii)$ 
$$\sum_{i=N}^{\infty} (a_{i+1}-a_{i}) i^k=O(2^{-N}).$$
This allows us to calculate the contribution
from $S_1$ up to $O(1/g^r)$. The other two terms $S_2$ and $S_3$ can be treated in a similar way.
\hfill $\Box$

In order to prove Theorem $\ref{theo:general:2}$ we need two auxiliary  lemmas:

\begin{lemm}\label{lem:general:1}
\begin{enumerate}[(i)]
\item Fix $k\; (0<k\leq n)$ and $d_{1},\ldots, d_{k}\in \mathbb{Z}_{\geq 0}$. Then for ${d}=(d_{1},\ldots, d_{k}, 0,\ldots,0)$ each 
term $e_{n,{d}}^{(s)}$ in the asymptotic expansion (\ref{e}) is a polynomial in $n$ of degree at most $s$.
\item Each $a_{n}^{(s)}$ and $b_{n}^{(s)}$ in (\ref{theo:n}) and (\ref{theo:g}) is a polynomial in $n$ of degree $s$. 
\end{enumerate}
\end{lemm}

\noindent
{\bf Proof}. 
The proof is again by induction on $s$. We prove a slightly stronger version of $(i)$:

\noindent
($i'$) For given $k$ and $s$, there exist polynomials $q_{j}(d_{1},\ldots,d_{k}),\; j=0,\ldots, s,$ such that the term $e_{n,d}^{(s)}$
has the form
$$e_{n,d}^{(s)}=\sum_{j=0}^{s} e_{d,j} n^{j}$$ 
with $|e_{d,j}| \leq q_{j}(d_{1},\ldots,d_{k}).$ 
In other words, the coefficients of $e_{n,d}^{(s)}$ considered as a polynomial in $n$ grow at most polynomially 
in $d_1,\ldots, d_k$. This would imply that 
$$\sum_{i=1}^{\infty} (a_{i+1}-a_{i}) e_{d(i),j} <\infty,$$
where $d(i)=(d_1+i,d_{2},\ldots,d_{k},0,\ldots,0)$.

Now, by $(\ref{Obs:I})$ and $(\ref{Obs:II})$ the statement ($i'$) for $s=r$ implies part $(ii)$ of the Lemma for $s=r$,
that is clear in view of $(\ref{a:n:s})$ and $(\ref{next})$. 
Moreover, part $(ii)$ for $s=r$ and the statement ($i'$) for $s<r$ imply ($i'$) for $s=r.$
This follows from $(\ref{eq:main})$ by analyzing the contributions from the terms $S_{1},$ $S_{2}$ and $S_{3}$ 
as in $(\ref{eq:aa})$, $(\ref{eq:bb})$ and $(\ref{eq:cc})$. \hfill $\Box$

\begin{lemm}\label{lem:general:2}
\begin{enumerate}[(i)]
\item Let $n$ and $k$ be fixed, and let $d=(d_{1},\ldots,d_{k},d_{k+1},\ldots,d_{n})$ with 
$d_{k+1},\ldots, d_{n}$ fixed. Then there exists a polynomial
$P_s\in {\mathbb R}[x_{1},\ldots,x_{k}]$ (depending on $d_{k+1},\ldots, d_{n}$) of degree $2s$ such that 
$e_{n,d}^{(s)}= P_s(d_{1},\ldots,d_{k})$ provided $d_{j} \geq 2s$ for $j=1,\ldots,k$.  The coefficient at each monomial 
$d_{1}^{\alpha_{1}}\cdots d_{k}^{\alpha_k}$ in $P_s$ is a linear rational
combination of 
$\pi^{2s-2\left[\frac{|\alpha|+1}{2}\right]},\ldots,\pi^{-2s},$ where $|\alpha|=\alpha_1+\ldots+\alpha_k.$  Moreover, 
for arbitrary $d_j$ the difference $e_{n,{d}}^{(s)}- P_s(d_{1},\ldots,d_{k})$ is a linear rational combinations of 
$\pi^{2s-2|d|},\ldots, \pi^{-2s}.$ 

\item  Each $a_{n}^{(s)}$ and $b_{n}^{(s)}$ in (\ref{theo:n}) and (\ref{theo:g}) is a rational polynomial of degree at most $s$ in ${\mathbb Q}[\pi^2,\pi^{-2}].$
\end{enumerate}
\end{lemm}

\noindent
{\bf Proof.} 
The proof is by induction in $s$ and utilizes the same techniques that we have already used before, 
so we only sketch it here. The statement of the lemma follows from the following claims:

\noindent
{\bf Claim $1$}: {\em Part (ii) for $s<r$ implies that the coefficient at $1/g^{r}$ in the expansion of $V_{g_{1},n_{1}+1}\cdot 
V_{g-g_{1},n-n_{1}+1}/V_{g,n}$ is a polynomial of degree at most $r$ in ${\mathbb Q}[\pi^2, \pi^{-2}]$. More precisely, when 
$g_{1}, n_{1}$ and $n_2$ are fixed
$$\frac{V_{g_{1},n_{1}+1}\cdot V_{g-g_{1},n-n_{1}+1}}{V_{g,n}}= \sum_{k=2g_{1}+n_1-1}^s \frac{c^{(k)}_{g_1,n_1}}{g^{k}} + 
O\left(\frac{1}{g^{s+1}}\right),$$
where  $c^{(k)}_{g_1,n_1}$ is a polynomial of degree at most $s$ in ${\mathbb Q}[\pi^2, \pi^{-2}]$.}

\medskip
\noindent
This is a simple consequence of Remarks $\ref{rneed}$(2), $\ref{rem3}$(2) and formula $(\ref{subsurface})$. 
\medskip

\noindent
{\bf Claim $2$}:  {\em Part (i) of the lemma for $s=r$ and part (ii) for $s<r$ imply part (ii) for $s=r$.} 

\medskip
\noindent
Note that by part $(i)$, each element of the infinite sum $(\ref{a:n:s})$ is a polynomial of degree at most $r$ in ${\mathbb Q}
[\pi^2, \pi^{-2}].$ On the other hand, when $l\geq r$ the coefficient $e_{n,l}^{(r)}=P_r(l)$ is a polynomial in $l$. Remark 
$\ref{rneed}$(1) and the properties of the coefficients of $P_(l)$ imply that $a_{n}^{(r)}$ is a rational polynomial of degree at 
most $r$ in ${\mathbb Q}[\pi^2,\pi^{-2}].$  The similar statement for $b_{n}^{(r)}$ follows from $(\ref{Obs:I})$, Claim ${\bf 1}$ 
and part (i) for $s=r$ (when $k=1$ and $d=(1,0,\ldots,0)$).
\medskip

\noindent
{\bf Claim $3$}:  {\em Part (ii) for $s<r$ and part (i) for $s<r$ imply part (i) for $s=r$.}

\medskip
\noindent
First, we use $({\bf Ib})$ to find the polynomial $P_r$. 
We put $P_r(d_{1},\ldots,d_{k})= e^{(r)}_{n+2,d}$ for $d=(0,0,d_{1},\ldots,d_{k},d_{k+1},\ldots,d_{n})$ 
with $d_{1},\ldots, d_{k}\geq 2r$. Note that the number $[\tau_{l},\tau_{d_1}\ldots\tau_{d_k}]_{g,k+1} \neq 0$ only when 
$l \leq 3g-2+k-|d|,\; |d|=d_1+\ldots+d_k.$ On the other hand, by Lemma $\ref{estimate:upperbound}$, the term
$$\left[\tau_{0}^{2} \tau_{l} \; \prod_{i\in I } \tau_{d_{i}}\right]_{g_{1},|I|+3}  
\cdot \quad\left[\tau_{0}^{2}\prod_{i\in J} \tau_{d_{i}}\right]_{g_{2},|J|+2}$$
in $({\bf Ib})$ is of order $O(1/g^{s+1})$ unless either $2g_1+|I|+1 <s $ or $2g_2+|J| < s$. Similarly the term
$$ \left[ \tau_{0}\tau_{l} \; \prod_{i\in I } \tau_{d_{i}}\right]_{g_{1},|I|+2} 
\cdot \quad\left[\tau_{0}^{3}\prod_{i\in J} \tau_{d_{i}}\right]_{g_{2},|J|+3}$$
in $({\bf Ib})$ is of order $O(1/g^{s+1})$ unless  $2g_1+|I| <s $ or $2g_2+|J|+1 < s$. For both terms we can explicitly 
calculate the expansions of the factors up to $O(1/g^{s+1})$ if we know their expansions up to $O(1/g^{s}).$ 

Then Lemma $\ref{estimate:upperbound}$, formula $({\bf Ib})$ and the induction hypothesis imply that 
$$ P_s(d_{1}+1,d_2,\ldots,d_{k})= P_s(d_{1},\ldots,d_{k},0,0)+ P'_s(d_{1},\ldots,d_{k})+P''_s(d_{1},\ldots, d_{k}),$$
where $P'_s$ is a polynomial when $d_{1}\geq s-1,d_2\geq s,\ldots, d_{k} \geq s,$ and $P''_s(d_{1},\ldots,d_{k})$ 
is nontrivial only if $d_{1},\ldots, d_{k}\leq s.$ The result follows from Lemma 
$\ref{lem:general:1}(i)$ and Claim ${\bf 1}.$ 
In the simplest case $k=1,\;d_{2}=\ldots=d_{n}=0$, and $d_1=d\geq 2s$, 
the relation $({\bf Ib})$ implies that
$$[\tau_{d}\tau_0^{n-1}]_{g,n}=[\tau_{s}\tau_{0}^{n+2(d-s)-1}]_{g-(d-s),n+2(d-s)}+ Q(d),$$
where $Q$ is a polynomial in $d$.

Next, we use $({\bf III})$ to prove the statement about the coefficients of $P_s$. We explicitly calculate the expansion of
$[\tau_{d_1}\ldots\tau_{d_n}]_{g,n}$
using (\ref{eq:main}) and Remark $\ref{rem3}$. We evaluate the contributions from the terms $S_{1}, S_{2}$ and $S_{3}$ and show that there exist  polynomials $Q^{(s)}_{i},\; i=1,2,3,$ such that:
\begin{enumerate}
\item When $d_{1},\ldots, d_{k}$ are large enough (i.e., $\geq 2s$), the coefficient $S^{(s)}_{i}$ at $1/g^{s}$ in $S_i$
is equal to $Q^{(s)}_{i}(d_{1},\ldots,d_{k});$
\item The coefficient at the monomial $d_{1}^{\alpha_{1}}\cdots d_{k}^{\alpha_k}$ in each $Q^{(s)}_{i}$ is a linear rational 
combination of $\pi^{2s-2[(|\alpha|+1)/2]},\ldots,\pi^{-2s},$ where $|\alpha|=\alpha_1+\ldots+\alpha_k;$ 
\item  For all $d_{1},\ldots, d_{k}$, the difference $ S_{i}^{(s)}- Q^{(s)}_{i}(d_{1},\ldots,d_{k})$ 
is a linear rational combinations of $\{\pi^{2s-2|d|},\ldots, \pi^{-2s}\}.$ 
\end{enumerate}

\noindent
{\em Contributions from $S_1$.}
We use the induction hypothesis and Lemma $\ref{aa}$ to expand the ratio $\widetilde{{A}}_{{d},g,n}/V_{g,n-1}$ up to the order $O(1/g^s)$ by expanding each  $[\tau_{d_{1}+d_{j}+i-1}\tau_{d_{2}}\ldots\widehat{\tau_{d_{j}}}\ldots \tau_{d_{n}}]_{g,n-1}/V_{g,n-1}$ up to $O(1/g^{s-1}).$
Put $$q^{(m)}(d) =8 \sum_{j=2}^{n} \,\sum_{i=0}^{\infty} (a_{i+1}-a_{i}) e_{n,{d_{j}(i)}}^{(m-1)},$$ 
where $d_{j}(i)=(d_{1}+d_{j}+i-1,d_{2},\ldots,, \widehat{d_j}\ldots, d_{n}).$
Then 
following Remark $\ref{rneed}$(2), and $(\ref{tf})$ 
$$S_{1}= \sum_{j=1}^{s} \frac{Q_{1}^{(j)}(d)}{g^{j}}+O\left(\frac{1}{g^{s}}\right),$$
where 
$$Q_{1}^{(j)}(d)=\sum_{j_{1}+j_{2}=j}  q^{(j_1)}(d) \cdot a_{n-1}^{(j_2)}.$$

Now by the induction hypothesis and part $(ii)$ of the lemma for $s<r$, both $q^{(m)}(d)$ and $Q_{1}^{(m)}(d)$
are polynomials in $d_1,\ldots, d_k$ of degree $2m$ whenever $d_{1},\ldots, d_{k}$ are large enough.
The coefficient at $d_{1}^{\alpha_{1}}\ldots d_{k}^{\alpha_{k}}$ in $q^{(m)}(d)$ is a rational linear combination of 
the terms of the form $c_{\alpha_j(r)} \sum_{i=0}^{\infty} i^{r} (a_{i+1}-a_{i}),$ where $c_{\alpha_j(r)}$ is the coefficient at 
$x_{1}^{r+\alpha_{1}+\alpha_{j}} x_{2}^{\alpha_{2}} \ldots\widehat{x_j^{\alpha_j}}\ldots x_{k}^{\alpha_k}$ in $P^{(m-1)}(x_{1},\ldots,\widehat{x_{j}},\ldots,x_{k})$. Now by Lemma $\ref{aa},\;(iv)$ the sum $\sum_{i=0}^{\infty} i^{r} (a_{i+1}-a_{i})$
is a rational linear combination of $\pi^{2}, \ldots, \pi^{2[r/2]}$, and by the induction hypothesis $c_{\alpha_j(r)}$ is a rational linear combination of $\pi^{2m-2-2r-|\alpha|} \ldots, \pi^{-2m-2}$. Therefore the coefficients of these polynomials 
are rational combinations of $\pi^{2m-|{\alpha}|},\ldots, \pi^{-2m}.$

\medskip
\noindent
{\em Contributions from $S_2$.}
In the same way, the induction hypothesis and Lemma $\ref{aa}$ allow to expand $\widetilde{{B}}_{d,g,n}/V_{g-1,n+1}$ 
up to the order of $1/g^s.$ Repeating the proof of case ({\bf 2}) of Theorem $\ref{theo:tauk}$ we see that each term 
in the expansion is of the form $(\ref{eq:bb})$.
In the expansion of $T_{m,g,n}$ the term of order $1/g^{m}$
is a polynomial in $d_1,\ldots, d_n$ of degree $m+1,$ whose coefficients satisfy the properties $1,$ $2$, and $3$ 
mentioned above.

\medskip
\noindent
{\em Contributions from $S_3$.}
What is different here compared to step {\bf 3} in the proof of Theorem $\ref{theo:tauk}$, is that 
$\widetilde{{C}}_{{d},g,n}$ contributes to the terms of order $1/g^{2},\ldots,1/g^{r}.$
However, by Lemma $\ref{estimate:upperbound}$ these contributions can be evaluated using $(\ref{subsurface}).$
More precisely, as in Remark $\ref{rem3}$, we can write $S_3$ as a sum of finitely many elements of the form
\begin{align*} 
\frac{\left[\tau_{k} \prod_{i\in I } \tau_{d_{i}}\right]_{g',|I|+1}}{V_{g', |I|+1}} 
\cdot &\frac{V_{g-g',|J|+1}\cdot V_{g', |I|+1}}{V_{g,n}}\times\\
\times &\sum_{i=0}^{\infty} (a_{i+1}-a_{i}) \;  \frac{\left[\tau_{d_1+i-2-k} \prod_{i\in J} \tau_{d_{i}}\right]_{g-g',|J|+1}}{V_{g-g',|J|+1}},
\end{align*}
where $I \amalg J=\{2,\ldots,n\},$ $0\leq g' \leq g,$ and $2g'-1+|I| \leq r.$ The result now follows from the induction
hypothesis on the behavior of the coefficient at $1/g^{s}$ in the expansion of
 $\frac{\left[\tau_{d_1+i-2-k} \prod_{i\in J} \tau_{d_{i}}\right]_{g-g',|J|+1}}{V_{g-g',|J|+1}}$ for $s<r$, together
Claim $1$ and Remark $\ref{rneed}(2),(3).$
\hfill $\Box$

\medskip
\noindent
{\bf Proof of Theorem $\ref{theo:general:2}$.}  
It is easy to see that part $(i)$ of the theorem is a special case of Lemma $\ref{lem:general:2}$ for $n=k.$ 
Part $(ii)$ is a consequence of Lemma $\ref{lem:general:2}(ii)$ and $\ref{lem:general:1}(ii).$
\hfill $\Box$

\medskip
\noindent
{\bf Proof of Theorem $\ref{theo:main:fixn}$.}
For a fixed $n \geq 0$, Theorem $\ref{theo:general:1} (ii)$ applied to the obvious identity $\frac{V_{g+1,n}}{V_{g,n}}=\frac{V_{g+1,n}}{V_{g,n+2}}\cdot\frac{V_{g,n+2}}{V_{g,n+1}}\cdot\frac{V_{g,n+1}}{V_{g,n}}$ immediately yields
\begin{align*} 
\frac{V_{g+1,n}}{V_{g,n}}= &(4\pi^2)^2(2g+n-1)(2g+n-2)\times\\ 
\times &\left(1-\frac{1}{2g}\right)\cdot\left(1+ \frac{r^{(2)}_{n}}{g^2}+  \ldots +\frac{r^{(s)}_{n}}{g^s}+ O\left(\frac{1}{g^{s+1}}\right)\right)
\end{align*}
as $g \rightarrow \infty.$
On the other hand, we have
$$ V_{g,1}= \prod_{j=2}^{g-1}  \frac{V_{j+1,1}}{V_{j,1}} \cdot V_{2,1},$$
and therefore the result of Theorem $\ref{theo:main:fixn}$ for $n=1$ is a consequence of the following:
\begin{lemm}
Let $a_{2},\ldots, a_{l} \in {\mathbb R}$ and let $\{c_{j}\}_{j=1}^{\infty}$ be a positive sequence with the property
$$ c_{j}=1+\frac{a_{2}}{j^{2}}+\ldots+\frac{a_{s}}{j^{s}}+ O\left(\frac{1}{j^{s+1}}\right).$$
Then there exist $b_{1},\ldots, b_{s-1}$ such that 
$$ \prod_{j=1}^{g} c_{j}= C_0 \; \left(1+\frac{b_{1}}{g}+\ldots
\frac{b_{s-1}}{g^{s-1}}+ O\left(\frac{1}{g^{s}}\right)\right), $$
as $g \rightarrow \infty,$ 
where $C_0=\prod_{j=1}^{\infty} c_{j}.$ Moreover, $b_{1},\ldots,b_{s-1}$ are polynomials  in $a_{2},\ldots, a_{s}$ with rational coefficients.
\end{lemm}

\noindent
{\bf Proof.}
First, we can write 
$$
c_{j}=R_{j}\left(1+\frac{d_{2}}{j^2}\right) \ldots \left(1+\frac{d_{l}}{j^l}\right),
$$
where $d_{2},\ldots d_{l}$ are polynomials in $a_{2},\ldots, a_{l}$ and  $R_{j}=1+ O(1/j^{l+1})$ as $j\to\infty$.
One can check that there exists a constant $R$ such that 
 $$\prod_{j=1}^{g} R_{j}=R \; \left(1+O\left(\frac{1}{g^l}\right)\right)$$
as $g\to\infty$. So it is enough to prove that
$$\prod_{j=1}^{g} \left(1+\frac{d_{k}}{j^{k}}\right)= C_{k} \; \left(1+\frac{p_{k}^{(1)}}{g}+\ldots
\frac{p_{k}^{(l-1)}}{g^{l-1}}+ O\left(\frac{1}{g^{l}}\right)\right) $$
as $g\to\infty$, where $p_{k}^{(1)},\ldots, p_{k}^{(l)}$ are polynomials in $d_{k}$ with rational coefficients, and $C_{k}= \prod_{j=1}^{\infty} \left(1+\frac{d_{k}}{j^{k}}\right),\; k\geq 2.$
It is enough to get bounds for the error term 
\begin{equation}\label{BB}
E_{N,k}=\log \prod_{j=N}^{\infty} \left(1+\frac{d_{k}}{j^{k}}\right)= \sum_{j=N}^{\infty} \log\left(1+ \frac{d_k}{j^k}\right).
\end{equation}
Using the Taylor series $\log(1+x)= x-x^{2}/2+x^{3}/3-\ldots$, we expand each term $\log\left(1+ \frac{d_k}{j^k}\right)$ up to the order $j^{-l}$. 

Now we make use of the Euler-Maclaurin summation formula (cf. \cite{E}):
\begin{align*}
\zeta(s)=&\sum_{i=1}^{N} \frac{1}{i^{s}}+\frac{1}{(s-1) N^{s-1}}+ \frac{1}{2N^s}+\\
+&\sum_{m=1}^{r-1} \frac{s(s+1)\ldots (s+2m-2)B_{2m}}{(2m)! \, N^{s+2m-1}}+E_{2r}(s),
\end{align*}
where $\zeta(s)$ is the Riemann zeta function, $B_{2m}=(-1)^{m+1} 2 (2m)!  \frac{\zeta(2m)}{(2 \pi)^{2m}}$
is the $m$th Bernoulli number, and the error term $E_{2r}(s)$ has an estimate
$$|E_{2r}(s)|<\left| \frac{s(s+1)\ldots (s+2r-2)B_{2r}}{({\rm Re}(s) +2r-1)(2r)! \, N^{s+2r-1}}\right|$$
(the formula holds for all $s$ with ${\rm Re}(s)>-2r+1$). 
From here we get that for any integer $k>1$
\begin{align*}
\sum_{i=N+1}^{\infty} \frac{1}{i^{k}}=&\frac{1}{(k-1) N^{k-1}}+\frac{1}{2N^k}+\\
+&\sum_{m=1}^{[(l-k)/2]} \frac{B_{2m}\cdot k(k+1)\ldots (k+2m-2)}{(2m)! \; N^{k+2m-1}}
+O\left(\frac{1}{N^l}\right)
\end{align*}
as $N\to \infty$.
Therefore, given $l$, there exist $q_{k}^{(k-1)},\ldots ,q_{k}^{(l)}$ such that 
$$E_{N,k}= \frac{q_{k}^{(k-1)}}{N^{k-1}}+\ldots+\frac{q_{k}^{(l)}}{N^l}+ O\left(\frac{1}{N^{l+1}}\right),$$
where each $q_{k}^{(i)}$ is a polynomial with rational coefficients in $d_{k}.$ Now we can easily get hold on the error terms $e^{E_{N,k}}$:
$$e^{E_{N,k}}= \frac{p_{k}^{(k-1)}}{N^{k-1}}+\ldots+\frac{p_{k}^{(l)}}{N^l}+ O\left(\frac{1}{N^{l+1}}\right),$$
where $p_{k}^{(j)}$ is a polynomial in $q_{k}^{(i)},\; i=k-1,\ldots,l$, which implies the result.
\hfill $\Box$
As a result, we can write 
 $$V_{g,1} =C\,\frac{(2g-2)!\,(4\pi^2)^{2g-2}}{\sqrt{g}} \left(1+\frac{c_1^{(1)}}{g}+\ldots+\frac{c_1^{(k)}}{g^{k}}+ O\left(\frac{1}{g^{k+1}}\right)\right),$$
 as $g \rightarrow \infty.$ On the other hand,
$$\frac{V_{g,n}}{C\cdot C_{g,n}}= \prod_{j=1}^{n-1}  \frac{V_{g,j+1}}{4\pi^2 (2g-2+j) V_{g,j}} \cdot \frac{V_{g,1}}{C \cdot C_{g,n}},$$
where $C_{g,n}=(2g-3+n)!\,(4\pi^2)^{2g-3+n}\,g^{-1/2},$ $C=\lim_{g\rightarrow \infty} \frac{V_{g,1}}{C_{g,1}}.$
The following is elementary:

\noindent
{\it Fact.}
Let $\{f_{i}\}_{i=1}^{\infty}$ be a sequence of functions with the expansion $$f_{i}(g)=1+\frac{p(1,i)}{g}+\ldots+\frac{p(s,i)}{g^s}+O(\frac{1}{g^{s+1}}).$$
Assume that for a given $j$, $p(j,k)$ is a polynomial in $k$ of degree $j$. 
Let $$H(g,n)=\prod_{j=1}^{n} f_{j}(g).$$ Then 
$$ H(g,n)=1+\frac{h_1(i)}{g}+\ldots+\frac{h_{s}(i)}{g^s}+O(\frac{1}{g^{s+1}}),$$
where for a given $j$, $h_{j}(k)$ is a polynomial in $k$ of degree $2j$. Moreover the leading coefficient of $h_{j}(k)$ is equal to $\frac{l^{j}}{2^j j!},$
where $l$ is the leading coefficient of the linear polynomial $p(1,j).$

This fact is a consequence of an elementary observation. Given $m_{1},\ldots m_{k} \in {\mathbb N}$, let $$F_{m_1,\ldots,m_k}(d)=\sum_{\sumind{x_{1},\ldots,x_{k} \in \{1,\ldots, n\},}  {{x_{i}\neq x_{j}}}} x_{1}^{m_{1}}\cdots x_{k}^{m_{k}}$$ is a polynomial in $d$ of degree $(m_{1}+1)+\ldots (m_{k}+1)$. Note that for $d<k,$ $F_{m_1,\ldots,m_k}(d)=0.$

Now we use this observation for $f_{j}(g)= \frac{V_{g,j+1}}{4\pi^2 (2g-2+j) V_{g,j}}.$ In this case, $(\ref{eq:n})$ implies that $l= -1/\pi^2$. Hence the result follows from $(\ref{theo:n})$ and Theorem $\ref{theo:general:2},\;(ii)$.
\end{section}

\begin{section}{Asymptotics for variable $n$}
In this section we discuss the asymptotics behavior of $V_{g,n(g)}$ in case when $n(g)\rightarrow \infty $ as $g \rightarrow \infty$ and prove Theorem $\ref{Main}:$ if 
$n(g)^2/g \rightarrow 0$ as $g \rightarrow \infty,$ then
$$\lim_{g\rightarrow \infty} \frac{V_{g,n(g)}}{C_{g,n(g)}}=C,$$
where $C_{g,n}=(2g-3+n)!\,(4\pi^2)^{2g-3+n}\,g^{-1/2}$ and $C=\lim_{g\rightarrow \infty} \frac{V_{g,0}}{C_{g,0}}$.

We need the following basic lemma:
\begin{lemm}\label{lemm:general:n}
There are universal constants $c_{0},c_{1},c_{2}, c_{3}, c_4>0$ such that for $g,n \geq 0$ the following inequalities hold:

\begin{enumerate}[(i)]
\item  for any $k \geq 1,$ $$ c_{0}\cdot \frac{n}{2g-2+n} \leq 1-\frac{[\tau_{k}\,\tau_{0}^{n-1}]_{g,n}}{V_{g,n}} \leq c_{1}\cdot \frac{n k^{2}}{2g-2+n},$$
\item $$ \left| \frac{(2g-2+n) V_{g,n}}{V_{g,n+1}}-\frac{1}{4\pi^2}\right|  \leq c_{2}\cdot \frac{n}{2g-2+n},$$ 
\item $$  \frac{V_{g-1,n+4}}{V_{g,n+2}} \leq 1-c_{4}\cdot\frac{n}{2g-2+n}.$$
\end{enumerate}
\end{lemm}
 \noindent
 {\bf Proof.} The proof follows the same lines as the proofs of formulas $(\ref{f1}), (\ref{f2})$ and $(\ref{simple})$ (see \cite{M:large} for details). 
 First, observe that 
 \begin{itemize} 
 \item Recursion ${(\bf III)}$ implies 
  \begin{equation}\label{simple:bound:1}
  [\tau_{k+1}\,\tau_{0}^{n-1}]_{g,n}\leq  [\tau_{k}\,\tau_{0}^{n-1}]_{g,n} \leq V_{g,n},
 \end{equation} 
and we also have  
$$b \leq \frac{[\tau_{1}\,\tau_{0}^{n-1}]_{g,n}}{ V_{g,n}}$$
 where $b=\max\{a_{i}/a_{i+1}\},\; i=0,1,\ldots;$
 \item For $l\geq 0$
 $$\frac{l\,\pi^{2l-2}}{(2l+1)!} \geq \frac{(l+1)\pi^{2l}}{(2l+3)!},$$
 so that from $({\bf II})$ and $(\ref{simple:bound:1})$ we have 
\begin{equation}\label{observation:II}
b_{0} \leq \frac{(2g-2+n) V_{g,n}}{V_{g,n+1}} \leq b_{1},
\end{equation}
where
$$b_{0}=b\cdot\left(\frac{1}{6}-\frac{\pi^2}{60}\right), \qquad 
b_{1}= \sum_{l=1}^{\infty}  \frac{l\,\pi^{2l-2}}{(2l+1)!}.$$
 \item By recursion $({\bf Ia})$, for $n \geq 2$ 
 \begin{equation}\label{observation:I}
\frac{V_{g-1,n+4}}{V_{g,n+2}} \leq 1.
 \end{equation}
 \end{itemize}
Note that in view of $(\ref{observation:II})$ and $(\ref{observation:I})$, Theorem $\ref{theo:tauk},\;(ii)$ implies the first inequality. In order to prove $(ii)$
we will use $(i)$ and $(\ref{Obs:II})$. As a result, we get 
$$ \left| \frac{(2g-2+n) V_{g,n}}{V_{g,n+1}}-\frac{1}{4\pi^2}\right| \leq  \sum_{l=1}^{\infty} \frac{ l^3 \, \pi^{2l-2}}{(2l+1)!}\;\cdot  \frac{n}{2g-2+n}.$$

The third inequality is a simple consequence of $({\bf Ia})$. This bound can be obtained from the lower bound in $(\ref{observation:II})$ and $(\ref{Obs:I})$.
\hfill $\Box$

\noindent
{\bf Proof of Theorem $\ref{Main}$.}
By Lemma \ref{lemm:general:n},
$$\frac{(2g-2+n(g)) V_{g,n(g)}}{V_{g,n(g)+1}} \rightarrow \frac{1}{4\pi^2},\qquad 
\frac{V_{g-1,n(g)+4}}{V_{g,n(g)+2}} \rightarrow 1,$$ 
when $n(g)/g \rightarrow 0$ as $g \rightarrow \infty.$
From definition of $C_{g,n}$ immediately follows that 
$$\frac{V_{g,n(g)}}{C_{g,n(g)}}=\frac{V_{g,n(g)}}{4\pi^{2} (2g-2+n(g)) V_{g,n(g)-1}}\cdot\frac{V_{g,n(g)-1}}{4\pi^{2}(2g-3+n(g))  V_{g,n(g)-2}}\cdots\frac{V_{g,0}}{C_{g,0}}. $$
Hence, by Theorem $\ref{theo:main:fixn}$ and Lemma $\ref{lemm:general:n},\;(ii)$ we have 
\begin{align*} 
(1-c_{2}/g)\cdot (1-2\cdot c_{2}/g)& \cdots (1-n(g)\cdot c_{2}/g)\cdot (C+O(1/g))\leq\\
\leq&\frac{V_{g,n(g)}}{C_{g,n(g)}} \leq \\ 
\leq(1+c_{2}/g)\cdot (1+2\cdot c_{2}/g)& \cdots (1+n(g)\cdot c_{2}/g) (C+O(1/g)),
\end{align*}
which implies the result.
\hfill $\Box$
\end{section}

\noindent
{\bf Acknowledgements.} 

The work of MM is partially supported by an NSF grant.
The work of PZ is partially supported by the RFBR grants 11-01-12092-OFI-M-2011 and 11-01-00677-a, and he
gratefully acknowledges the hospitality and support of MPIM (Bonn), QGM (Aarhus) and SCGP (Stony Brook).
We thank Maxim Kazarian and Don Zagier for enlightening discussions.

\end{document}